\documentclass{amsart}

\newtheorem{theorem}{Theorem}[section]

\theoremstyle{definition}
\newtheorem{definition}[theorem]{Definition}

\theoremstyle{remark}

\usepackage{graphicx}
\usepackage{color}
\usepackage{amsmath}
\usepackage{amsfonts}
\usepackage{amssymb}
\newcommand{\be}{\begin{equation}}
\newcommand{\ee}{\end{equation}}

\newcommand{\al}{\alpha}

\newcommand{\bet}{\beta}

\newcommand{\om}{\omega}

\newcommand{\bma}{\begin{pmatrix}}
\newcommand{\ema}{\end{pmatrix}}

\newcommand{\dz}{\wedge}

\newcommand{\ba}{\begin{array}}

\newcommand{\ea}{\end{array}}

\newcommand{\beq}{\begin{eqnarray}}

\newcommand{\eeq}{\end{eqnarray}}

\newtheorem{lm}{lemma}

\newtheorem{thee}{theorem}

\newtheorem{proo}{proposition}

\newtheorem{co}{corollary}

\newtheorem{rem}{remark}

\newtheorem{deff}{definition}

\newcommand{\bd}{\begin{deff}}

\newcommand{\ed}{\end{deff}}

\newcommand{\bl}{\begin{lm}}

\newcommand{\el}{\end{lm}}

\newcommand{\bp}{\begin{proo}}

\newcommand{\ep}{\end{proo}}

\newcommand{\bt}{\begin{thee}}

\newcommand{\et}{\end{thee}}

\newcommand{\bc}{\begin{co}}

\newcommand{\ec}{\end{co}}

\newcommand{\brm}{\begin{rem}}

\newcommand{\erm}{\end{rem}}

\newcommand{\der}{{\rm d}}

\usepackage{t1enc}
\def\frak{\mathfrak}

\newcommand{\newc}{\newcommand}

\let\ccdot\cdot
\def\cdot{\hbox to 2.5pt{\hss$\ccdot$\hss}}

\newc{\aR}{\mbox{\boldmath{$ R$}}}
\newc{\aS}{\mbox{\boldmath{$ S$}}}
\newc{\aT}{\mbox{\boldmath{$ T$}}}
\newc{\aW}{\mbox{\boldmath{$ W$}}}

\newc{\aK}{\mbox{\boldmath{$ K$}}}
\newc{\aL}{\mbox{\boldmath{$ L$}}}

\usepackage{amssymb}
\usepackage{amscd}

\newcommand{\Rho}{{\mbox{\sf P}}}

\newc{\obstrn}[2]{B^{#1}_{#2}}

\newcommand{\rpl}                         
{\mbox{$
\begin{picture}(12.7,8)(-.5,-1)
\put(0,0.2){$+$}
\put(4.2,2.8){\oval(8,8)[r]}
\end{picture}$}}

\newcommand{\lpl}                         
{\mbox{$
\begin{picture}(12.7,8)(-.5,-1)
\put(2,0.2){$+$}
\put(6.2,2.8){\oval(8,8)[l]}
\end{picture}$}}

\usepackage{ifthen}

\newcommand{\bbR}{\mathbb{R}}

\newcommand{\soa}{\frak{so}}

\newcommand{\ga}{\gamma}

\newc{\tensor}[1]{#1}
\newc{\Mvariable}[1]{\mbox{#1}}
\newc{\down}[1]{{}_{#1}}
\newc{\up}[1]{{}^{#1}}

\newc{\JulyStrut}{\rule{0mm}{6mm}}
\newc{\midtenPan}{\mbox{\sf S}}
\newc{\midten}{\mbox{\sf T}}
\newc{\midtenEi}{\mbox{\sf U}}
\newc{\ATen}{\mbox{\sf E}}
\newc{\BTen}{\mbox{\sf F}}
\newc{\CTen}{\mbox{\sf G}}

\def\sideremark#1{\ifvmode\leavevmode\fi\vadjust{\vbox to0pt{\vss
 \hbox to 0pt{\hskip\hsize\hskip1em
 \vbox{\hsize3cm\tiny\raggedright\pretolerance10000
 \noindent #1\hfill}\hss}\vbox to8pt{\vfil}\vss}}}%

\numberwithin{equation}{section}

\newcounter{romenumi}
\newcommand{\labelromenumi}{(\roman{romenumi})}

\begin{document}

\title{Conformal structures with explicit ambient metrics and
  conformal $G_2$ holonomy} 

\author{Pawe\l~ Nurowski} \address{Instytut Fizyki Teoretycznej,
Uniwersytet Warszawski, ul. Ho\.za 69, Warszawa, Poland}\thanks{This
  work was supported in part by the Polish Ministerstwo Nauki i
  Informatyzacji grant nr: 1 P03B 07529 and the US Institute for Mathematics and
  Its Applications in Minneapolis}
\email{nurowski@fuw.edu.pl}

\date{\today}
\begin{abstract}
Given a generic $2$-plane field on a $5$-dimensional manifold we
consider its $(3,2)$-signature conformal metric $[g]$ as defined in
\cite{conformal}. Every conformal class $[g]$ obtained in this way has 
very special conformal holonomy: it must be contained in the
split-real-form of the exceptional group 
$G_2$. In this note we show that for special $2$-plane
fields on $5$-manifolds the conformal classes $[g]$
have the Fefferman-Graham ambient metrics which, contrary to the
general Fefferman-Graham metrics given as a formal power series
\cite{feff}, can
be written in an explicit form. We propose to study the 
relations between the conformal $G_2$-holonomy of metrics $[g]$ and
the possible pseudo-Riemannian $G_2$-holonomy of the corresponding 
ambient metrics. 
\end{abstract}
\maketitle

\allowdisplaybreaks

\section{The $(3,2)$-signature conformal metrics}
\noindent
Consider an equation
\be
z'=F(x,y,y',y'',z)\quad\quad{\rm  
with}\quad\quad F_{y''y''}\neq 0,\label{2m}\ee
for two real functions $y=y(x)$, $z=z(x)$ of one real variable $x$. To simplify notation introduce new symbols $p=y'$ and $q=y''$. Equation (\ref{2m}) is 
totally encoded in the system of three 1-forms:
\beq
\om^1&=&\der z-F(x,y,p,q,z)\der x\nonumber\\
\om^2&=&\der y-p\der x\label{fhil}\\
\om^3&=&\der p-q\der x,\nonumber
\eeq
living on a 5-dimensional manifold $J$ parametrized by $(x,y,p,q,z)$. In particular, every solution to (\ref{2m}) is a curve $\gamma(t)=(x(t),y(t),p(t),q(t),z(t))\subset J$ on which all the forms $\om^1,\om^2,\om^3$ identically vanish.

We introduce an equivalence relation between equations (\ref{2m}) which identifies the equations having the same set of solutions. This leads to the following definition:

\begin{definition} \label{equiv}
Two equations $z'=F(x,y,y',y'',z)$ and $\bar{z}'=\bar{F}(\bar{x},\bar{y},\bar{y}',\bar{y}'',\bar{z})$, defined on spaces $J$ and $\bar{J}$ parametrized, respectively, by $(x,y,p=y', q=y'',z)$ and $(\bar{x},\bar{y},\bar{p}=\bar{y}', \bar{q}=\bar{y}'',\bar{z})$, are said to be \emph{(locally) equivalent}, iff 
there exists a (local) diffeomorphism $\phi:J\to\bar{J}$ transforming the corresponding forms\\

\begin{tabular}{lcl}
$\om^1=\der z-F(x,y,p,q,z)\der x$&&$\bar{\om}^1=\der \bar{z}-\bar{F}(\bar{x},\bar{y},\bar{p},\bar{q},\bar{z})\der \bar{x}$\\
$\om^2=\der y-p\der x$&and\hspace{1cm}&$\bar{\om}^2=\der \bar{y}-\bar{p}\der \bar{x}$\\
$\om^3=\der p-q\der x$&&$\bar{\om}^3=\der \bar{p}-\bar{q}\der \bar{x}$
\end{tabular}\\

\noindent
via:\\

\begin{tabular}{c}
$\phi^*(\bar{\om}^1)~=~\al\om^1+\bet\om^2+\ga\om^3$\\
$\phi^*(\bar{\om}^2)~=~\delta\om^1+\epsilon\om^2+\lambda\om^3$,\\
$\phi^*(\bar{\om}^3)~=~\kappa\om^1+\mu\om^2+\nu\om^3$
\end{tabular} 
with functions $\al,\bet,\ga,\delta,\epsilon,\lambda,\kappa,\nu$ on $J$ such that 
$$
{\rm det}\bma
\al&\bet&\ga\\
\delta&\epsilon&\lambda\\
\kappa&\mu&\nu\\
\ema\neq 0.
$$
\end{definition}
It follows that equation (\ref{2m}) considered modulo equivalence relation of Definition \ref{equiv} uniquely defines a conformal class of $(3,2)$-signature metrics $[g_F]$ on the space $J$. In coordinates $(x,y,p,q,z)$ this class may be described as follows.
Let 
$$
D=\partial_x+p\partial_y+q\partial_p+F\partial_z
$$ 
be a total differential associated with equation (\ref{2m}) on $J$. Then a representative $g_F$ of the conformal class $[g_F]$ may be written as
\begin{eqnarray}
&&g_F=[~DF_{qq}^2 F_{qq}^2 + 6DF_q DF_{qqq}F_{qq}^2 - 6DF_{qqq}F_p F_{qq}^2
  -\nonumber\\&&3DDF_{qq} F_{qq}^3 + 9 DF_{qp}F_{qq}^3 - 9 F_{pp}F_{qq}^3 +\nonumber\\&&9
  DF_{qz} F_q F_{qq}^3 - 18 F_{pz} F_qF_{qq}^3 + 3DF_z
  F_{qq}^4 -\nonumber\\&& 6 DF_qF_{qq}^2 F_{qqp} + 6 F_p F_{qq}^2 F_{qqp} - 8
  DF_qDF_{qq}F_{qq}F_{qqq} + \nonumber\\&&8 DF_{qq}F_p F_{qq}F_{qqq} +
  3DDF_qF_{qq}^2 F_{qqq} - 3DF_p F_{qq}^2 F_{qqq} -\nonumber\\&& 3DF_z F_q F_{qq}^2
  F_{qqq} + 4(DF_q)^2 F_{qqq}^2 - 8 DF_q F_pF_{qqq}^2  -\nonumber\\&& 3 (DF_q)^2
  F_{qq}F_{qqqq}+ 4 F_p^2 F_{qqq}^2+ 6DF_q F_p
  F_{qq}F_{qqqq} -\nonumber\\&&  3 F_p^2 F_{qq}F_{qqqq} - 6 DF_q F_q F_{qq}^2
  F_{qqz} +6 F_p F_q F_{qq}^2 F_{qqz} - \nonumber\\&&3 DF_q F_{qq}^3 F_{qz}
  + 12 F_p F_{qq}^3 F_{qz} + 3 F_{qq}^2 F_{qqq}F_y -\nonumber\\&& 6
  DF_{qqq}F_qF_{qq}^2 F_z + 4DF_{qq} F_{qq}^3 F_z  + 6 F_q F_{qq}^2
  F_{qqp}F_z +\nonumber\\&& 8 DF_{qq}F_qF_{qq}F_{qqq}F_z - 4 DF_q
  F_{qq}^2 F_{qqq} F_z- \label{metr}\\&&9 F_{qp} F_{qq}^3 F_z + F_p F_{qq}^2 F_{qqq}
  F_z - 8 DF_q F_qF_{qqq}^2 F_z +\nonumber\\&& 8 F_p F_q F_{qqq}^2 F_z
  +
 6 DF_q F_q
  F_{qq}F_{qqqq} F_z - 6 F_p F_q F_{qq} F_{qqqq} F_z + \nonumber\\&&18
  F_{qq}^3 F_{qy} + 6 F_q^2 F_{qq}^2 F_{qqz} F_z + 3F_q F_{qq}^3
  F_{qz}F_z -\nonumber\\&& 2 F_{qq}^4 F_z^2 + F_q F_{qq}^2 F_{qqq} F_z^2 + 4 F_q^2
  F_{qqq}^2 F_z^2-\nonumber\\&&3 F_q^2 F_{qq} F_{qqqq}F_z^2 - 9 F_q^2
  F_{qq}^3 F_{zz}~]~(\tilde{\om}^1)^2+\nonumber\\&&\nonumber\\
&&[~6 DF_{qqq} F_{qq}^2 - 6 F_{qq}^2 F_{qqp} - 8
  DF_{qq}F_{qq}F_{qqq} + \nonumber\\&&
8 DF_q F_{qqq}^2 - 8 F_p F_{qqq}^2 - 6 DF_q F_{qq} F_{qqqq} + 
\nonumber\\&&6 F_p F_{qq}F_{qqqq} - 6 F_q F_{qq}^2F_{qqz} + 6 F_{qq}^3 F_{qz} +\nonumber\\&& 2
F_{qq}^2 F_{qqq} F_z - 8 F_q F_{qqq}^2 F_z + 6 F_q F_{qq}F_{qqqq}F_z~]~\tilde{\om}^1\tilde{\om}^2 +\nonumber\\&&\nonumber\\&& [~10 DF_{qq}F_{qq}^3 - 10 DF_q F_{qq}^2 F_{qqq} + 10 F_p F_{qq}^2 F_{qqq} -\nonumber\\&& 10 F_{qq}^4 F_z + 10 F_q F_{qq}^2 F_{qqq}F_z~]~\tilde{\om}^1\tilde{\om}^3 +\nonumber\\&&30 F_{qq}^4~\tilde{\om}^1\tilde{\om}^4 + [~30 DF_q F_{qq}^3 - 30 F_p F_{qq}^3 - 30 F_q F_{qq}^3 F_z~]~\tilde{\om}^1\tilde{\om}^5+\nonumber\\&&[~4 F_{qqq}^2 - 3 F_{qq} F_{qqqq}~]~(\tilde{\om}^2)^2-10 F_{qq}^2 F_{qqq}~\tilde{\om}^2\tilde{\om}^3 + 30 F_{qq}^3~\tilde{\om}^2\tilde{\om}^5 - 20 F_{qq}^4~(\tilde{\om}^3)^2\nonumber 
\end{eqnarray}
where\footnote{Note that formula for $g_F$ differs from the one given in Ref. \cite{conformal} by tilde signs over the all omegas. In Ref. \cite{conformal}, when copying the calculated metric $g_F$, by mistake, we forgot to put these tilde signs over the omegas. Hence, in Ref. \cite{conformal}, formula for $g_F$ is true, provided that one puts the tilde signs over the omegas and supplements it by the definitions (\ref{po}) of the tilded omegas.}

\begin{eqnarray}
\tilde{\om}^1&=&\der y-p\der x\nonumber\\
\tilde{\om}^2&=&\der z-F\der x-F_q(\der p-q\der x)\nonumber\\
\tilde{\om}^3&=&\der p-q\der x\label{po}\\
\tilde{\om}^4&=&\der q\nonumber\\
\tilde{\om}^5&=&\der x\nonumber.
\end{eqnarray}

It follows from the construction described in Ref. \cite{conformal} that when the equation (\ref{2m}) undergoes a diffeomorphism $\phi$ of Definition \ref{equiv}, the above metric $g_F$ transforms conformally.

The conformal class of metrics $[g_F]$ is very special among all the $(3,2)$-signature conformal metrics in dimension 5: the Cartan normal conformal connection for this class, instead of having values in full $\soa(4,3)$ Lie algebra, has values in its certain 14-dimensional subalgebra. This subalgebra turns out to be isomorphic to the split real form of the exceptional Lie algebra ${\mathfrak g}_2\subset\soa(4,3)$. Thus, conformal metrics $[g_F]$ provide an abundance of examples of metrics with an \emph{exceptional} conformal \emph{holonomy}. This holonomy is always a subgroup of the noncompact form of the exceptional Lie group $G_2$. We strongly believe that randomly chosen function $F$, such that $F_{qq}\neq 0$, give rise to conformal metrics $[g_F]$  with conformal holonomy \emph{equal} to $G_2$. 

It is interesting to study the conformal classes $[g_F]$ from the
point of view of the Fefferman-Graham ambient metric
construction \cite{feff}. Since for each $F$ defining equation (\ref{2m}) we have
a conformal class of metrics $[g_F]$ in dimension five, then since
five is \emph{odd}, Fefferman-Graham guarantees \cite{feff} that there is a
\emph{unique} formal power series of a \emph{Ricci-flat metric} of
signature $(4,3)$ corresponding to $[g_F]$. Moreover, since given $F$
the metric $g_F$ is explicitely determined by formula
(\ref{metr}), we see that starting with \emph{real analytic} $F$, the
metric $g_F$ is \emph{real analytic}. Thus, every analytic
$F$ of (\ref{2m}) leads to analytic $g_F$ and then, in turn, via
Fefferman-Graham, leads to a unique \emph{real analytic} ambient
metric $\tilde{g}_F$ of signature $(4,3)$. 
Since both 
the Levi-Civita connection for $\tilde{g}_F$ and the Cartan normal
conformal connection for the corresponding 5-dimensional metric $g_F$ 
have values in (possibly subalgebras of) the same Lie algebra
$\soa(4,3)$, it is
interesting to ask about the relations between them. We discuss these
relations on examples.   
\section{The strategy for constructing explicit examples of ambient
  metrics}\label{stra}
We start with the Fefferman-Graham result \cite{feff} adapted to the 5-dimensional
situation of conformal metrics $[g_F]$.

Let $g_F$ be a representative of the conformal class $[g_F]$ defined
on $J$ by (\ref{metr}). Consider a manifold
$J\times\bbR_+\times\bbR$. Introduce coordinates $(0<t,u)$ on 
$\bbR_+\times\bbR$ in $J\times\bbR_+\times\bbR$. 
We have a natural projection $\pi: J\times\bbR_+\times\bbR\to
J$, which enables us to pullback forms from $J$ to
$J\times\bbR_+\times\bbR$. Ommiting the pulback sign in the
expressions like $\pi^*(g_F)$ we define a formal power series
\be\label{ffg}
\check{g}_F=-2\der t\der u+t^2 g_F-ut\al+u^2\bet+u^3t^{-1}\gamma+\sum_{k=4}^\infty u^k t^{2-k}\mu_k. 
\ee  
Here $\al,\bet,\gamma,\mu_k$, $k=4,5,6,....$, are pullbacks of
symmetric bilinear forms $\al,\bet,\gamma,\mu_k$ 
from $J$ to $J\times\bbR_+\times\bbR$. Thus $\check{g}_F$ is a formal \emph{bilinear form} on
$J\times\bbR_+\times\bbR$. This formal bilinear form
has signature $(4,3)$ in some neighbourhood of $u=0$. The 
following theorem is due to Fefferman and
Graham \cite{feff}.
\begin{theorem}\label{fgt}
Among all the bilinear forms $\check{g}_F$ which, via (\ref{ffg}), are
associated with metric $g_F$ of (\ref{metr}) there is
\emph{precisely one}, say $\tilde{g}_F$, satisfying the Ricci flatness
condition $Ric(\tilde{g}_F)=0.$  
\end{theorem}
Given $g_F$, all the bilinear forms $\al,\beta,\gamma,\mu_k$ in
$\tilde{g}_F$ are totally
determined. Another issue is to calculate them explicitely. For
example, it is quite difficult to find the general formulas for the 
higher order forms $\mu_k$.  
Nevertherless the explicit expressions for the forms
$\al,\beta,\gamma$ are known \cite{grah,de}. We write them below in the
form obtained by C R Graham. We define the coefficients 
$\alpha_{ij}$, $\beta_{ij}$ and $\gamma_{ij}$ by  
$\alpha=\alpha_{ij}\der x^i\der x^j$, $\beta=\beta_{ij}\der x^i\der
x^j$, $\gamma=\gamma_{ij}\der x^i\der x^j$, 
where $(x^i)=(x,y,p,q,z)$ are coordinates on $J$. Then Graham's 
expressions for $\alpha_{ij}$, $\beta_{ij}$ and $\gamma_{ij}$ are \cite{grah}:
\begin{eqnarray}
&&\alpha_{ij}=2\Rho_{ij},\nonumber\\
&&\beta_{ij}=-B_{ij}+\Rho_i^{~k}\Rho_{jk},\nonumber\\
&&3\gamma_{ij}=B_{ij;k}^{~\quad k}-2W_{kijl}B^{kl}+4\Rho_{k(i}B_{j)}^{~~k}-4\Rho_k^{~k}B_{ij}+4\Rho^{kl}C_{(ij)k;l}-\label{ffgc}\\
&&2C^k_{~i}\phantom{}^lC_{ljk}+C_i^{~kl}C_{jkl}+2\Rho^k_{~k;l}C_{(ij)}^{\quad
    l}-2W_{kijl}\Rho^k_{~m}\Rho^{ml},\nonumber
\end{eqnarray}
where
$$\Rho_{ij}=\tfrac{1}{3}(R_{ij}-\tfrac{1}{8}Rg_{Fij}),$$
is the Schouten tensor for the metric $g_F=g_{Fij}\der x^i\der x^j$,
$$W_{ijkl}=R_{ijkl}-2(\Rho_{i[k}g_{Fl]j}-\Rho_{j[k}g_{Fl]i})$$
is its Weyl tensor,
$$C_{ijk}=\Rho_{ij;k}-\Rho_{ik;j}$$
is the Cotton tensor, and 
$$B_{ij}=C_{ijk;}^{\quad k}-\Rho^{kl}W_{kijl}$$
is the Bach tensor. 

Of course all the above quantities can be
explicitely calculated once $F$, and in turn the metric $g_F$, is chosen.

In the rest of the paper we will chose particular functions 
$F=F(x,y,p,q,z)$, and we will
calculate the corresponding forms $\alpha,\beta,\gamma$ for them. We
will give
examples of $F$'s for which the bilinear form $\gamma$ is identically
vanishing, \be
\gamma\equiv 0.\label{gam}
\ee 
Given such $F$'s we will consider  
$$\bar{g}_F=-2\der t\der u+t^2 g_F-ut\al+u^2\bet.$$
Note that $\bar{g}_F$ coincides with the ambient metric
$\tilde{g}_F$ up to the terms \emph{quadratic} in the ambient
coordinates $t,u$. If by \emph{chance} the bilinear form $\bar{g}_F$ 
satisfies the Ricci flatness condition 
$$Ric(\bar{g}_F)\equiv 0,$$ then by the
\emph{uniqueness} of the ambient metric $\tilde{g}_F$ stated in
Theorem \ref{fgt}, it will \emph{coicide} with the ambient metric
$\tilde{g}_F$:
$$\bar{g}_F\equiv \tilde{g}_F.$$ The uniqueness result of Theorem \ref{fgt},
together with the Ricci flatness of $\bar{g}_F$, is powerfull enough to
guarantee that not only the coefficient $\gamma$ in the ambient metric
$\tilde{g}_F$ identically vanishes, but that \emph{all} the
coefficients $\mu_k$, $k=4,5,6,....,$ vanish too!

Thus the strategy of finding explicit ambient metrics $\tilde{g}_F$
for $g_F$ is as follows:
\begin{itemize}
\item find $F=F(x,y,p,q,z)$ for which the corresponding metric $g_F$
  has identically vanishing form $\gamma$ of (\ref{ffgc});
\item calculate the approximate ambient metric $\bar{g}_F$ for such $F$;
\item check if the Ricci
  tensor $Ric(\bar{g}_F)$ of $\bar{g}_F$ is identically vanishing;
\item if you have $F$ with the above properties then the approximate metric
  $\bar{g}_F$ is the ambient metric $\tilde{g}_F$ for $g_F$.
\end{itemize} 

\section{Conformally Einstein example}
As the first example, following Ref. \cite{conformal}, we calculate
$g_F$ and its approximate ambient metric $\bar{g}_F$ for a very simple 
equation:
$$z'=F(y''),\quad\quad{\rm with}\quad\quad F_{y''y''}\neq 0.$$ 
It was shown in Ref. \cite{conformal} that the conformal class $[g_F]$
may be represented by\footnote{The metric presented here differs from
  this of \cite{conformal} by a convenient conformal factor equal to
  $-15(F'')^{10/3}$.} 
\beq
&&-15(F'')^{10/3}g_F=\nonumber\\
&& 30(F'')^4~[~\der q\der y-p\der q\der
  x~]~+~[~4F^{(3)2}-3F''F^{(4)}~]~\der z^2+\nonumber\\&&
2~[-5(F'')^2F^{(3)}-4F'F^{(3)2}+3F'F''F^{(4)}~]~\der p\der z+\nonumber\\&&
2~[15(F'')^3+5q(F'')^2F^{(3)}-4FF^{(3)2}+4qF'F^{(3)2}+3FF''F^{(4)}-\nonumber\\&&3qF'F''F^{(4)}~]~\der
  x\der z+\nonumber\\&&
[-20(F'')^4+10F'(F'')^2F^{(3)}+4(F')^2F^{(3)2}-3(F')^2F''F^{(4)}~]~\der
  p^2+\label{metprzy}\\&&
2~[-15F'(F'')^3+20q(F'')^4+5F(F'')^2F^{(3)}-10qF'(F'')^2F^{(3)}+\nonumber\\&&
4FF'F^{(3)2}-4q(F')^2F^{(3)2}-3FF'F''F^{(4)}+3q(F')^2F''F^{(4)}~]~\der
  p\der x+\nonumber\\&&
[-30F(F'')^3+30qF'(F'')^3-20q^2(F'')^4-\nonumber\\&&10qF(F'')^2F^{(3)}+10q^2
  F'(F'')^2F^{(3)}+4F^2F^{(3)2}-\nonumber\\&&
8qFF'F^{(3)2}+4q^2(F')^2F^{(3)2}-3F^2F''F^{(4)}+\nonumber\\&&6qFF'F''F^{(4)}-3q^2(F')^2F''F^{(4)}~]~\der
x^2.\nonumber
\eeq
As noted in Ref. \cite{conformal} this metric is conformal to a Ricci flat metric $\hat{g}_F={\rm e}^{2\Upsilon(q)}g_F$ with a conformal scale $\Upsilon=\Upsilon(q)$ satisfying second order ODE:
$$90F''^2(\Upsilon''-\Upsilon'^2)-60F''F^{(3)}\Upsilon'+3F''F^{(4)}-4F^{(3)2}=0.$$
Thus, since for each $F=F(q)$ the conformal class $[g_F]$ contains a Ricci flat metric, its conformal holonomy must be a proper subgroup of the noncompact form of $G_2$. 
An interesting feature of this conformal class is that it is very special among all the conformal classes associated with equation (\ref{2m}). Not only has $g_F$ very special conformal holonomy, making it very similar to the Lorentzian 4-dimensional Brinkman  metrics; moreover, since its Weyl tensor has essentially only one nonvanishing component (see Ref. \cite{conformal} for details) it is \emph{not} weakly generic (see Ref. \cite{gover} for definition). This makes $[g_F]$ analogous to the Lorentzian type $N$ metrics in 4-dimensions, such as for example, Fefferman metrics. 

Having $g_F$ of (\ref{metprzy}) we used the symbolic computer
calculation program Mathematica to
calculate its associated form $\gamma$ of (\ref{ffgc}). 
We checked that this form \emph{identically vanishes}. We further used
Mathematica to calculate the corresponding approximate ambient metric
$\bar{g}_F$. On doing that we obsereved that, surprisingly, the
bilinear form $\beta$ is also \emph{identically vanishing}. The
explicit formula for the approximate ambient metric is given below:
\be
\bar{g}_F=t^2 g_F-2~\der t\der u~-~2tuF''^{4/3}P\der q^2,\label{mmep}
\ee
with 
$$P=\frac{4F^{(3)2}-3F''F^{(4)}}{90(F'')^{10/3}},$$
and $g_F$ given by (\ref{metprzy}).
The metric $\bar{g}_F$ is defined locally on $J\times\bbR_+\times\bbR$
with coordiantes $(x,y,p,q,z,t,u)$. It obviously has signature $(4,3)$.  
We also checked, again using Mathematica, that $Ric(\bar{g}_F)\equiv 0$.
Thus, we fulfiled the strategy outlined in Section \ref{stra}. This
enables us to 
conclude that $\bar{g}_F$ of (\ref{mmep}) 
coincides with the ambient metric $\tilde{g}_F$ for $g_F$. 
To give expressions for the Cartan normal conformal connection for
$g_F$ and the Levi-Civita connection for $\tilde{g}_F=\bar{g}_F$ we first
introduce a nonholonomic coframe 
$(\theta^1,\theta^2,\theta^3,\theta^4,\theta^5)$ on $J$ given by
\begin{eqnarray*}
&&\theta^1=\der y -p\der x\\
&&\theta^2=\der z-F\der x-F'(\der p-q\der x)\\
&&\theta^3=-\tfrac{2}{\sqrt{3}}(F'')^{1/3}(\der p-q\der x)\\
&&30
(F'')^{10/3}\theta^4=\big(3F'F''F^{(4)}-4F'F^{(3)2}-10(F'')^2F^{(3)}\big)\big(\der
p-q\der x\big)+\\
&&\big(4F^{(3)2}-3F''F^{(4)}\big)\big(\der z-F\der x\big)+30(F'')^3\der x\\
&&\theta^5=-(F'')^{2/3}\der q.
\end{eqnarray*}
In this coframe the metric $g_F$ is simply:
$$g_F=2\theta^1\theta^5-2\theta^2\theta^4+(\theta^3)^2.$$
By means of the canonical projection
$$\pi(x,y,p,q,z,t,u)=(x,y,p,q,z)$$
the coframe $(\theta^1,\theta^2,\theta^3,\theta^4,\theta^5)$ can be pulbacked to five linearly independent forms
$(\theta^1,\theta^2,\theta^3,\theta^4,\theta^5)$ on
$J\times\bbR_+\times\bbR$. They can be suplemented by $$\theta^0=\der
t\quad\quad{\rm and}\quad\quad \theta^6=\der u$$ to form a coframe
$(\theta^0,\theta^1,\theta^2,\theta^3,\theta^4,\theta^5,\theta^6)$ on
the ambient space $J\times\bbR_+\times\bbR$.

The Cartan normal conformal connection, when written on $J$ in the coframe
$(\theta^1,\theta^2,\theta^3,\theta^4,\theta^5)$ reads:
$$
\om_{G_2}=
\bma
0&0&0&0&0&-P\theta^5&0\\
&&&&&&\\
\theta^1&0&Q\theta^2+\tfrac{9}{2\sqrt{3}}P\theta^3&\tfrac{1}{\sqrt{3}}\theta^4&-\tfrac{1}{2\sqrt{3}}\theta^3&0&-P\theta^5\\&&&&&&\\
\theta^2&0&0&\tfrac{1}{\sqrt{3}}\theta^5&0&-\tfrac{1}{2\sqrt{3}}\theta^3&0\\&&&&&&\\
\theta^3&0&-2\sqrt{3}P\theta^5&0&\tfrac{1}{\sqrt{3}}\theta^5&-\tfrac{1}{\sqrt{3}}\theta^4&0\\&&&&&&\\
\theta^4&0&0&-2\sqrt{3}P\theta^5&0&Q\theta^2+\tfrac{9}{2\sqrt{3}}P\theta^3&0\\&&&&&&\\
\theta^5&0&0&0&0&0&0\\&&&&&&\\0&\theta^5&-\theta^4&\theta^3&-\theta^2&\theta^1&0
\ema.
$$
Here:
$$Q=\frac{40F^{(3)3}-45F''F^{(3)}F^{(4)}+9F''^2F^{(5)}}{90F''^5}.$$
Now we use coframe
$(\theta^0,\theta^1,\theta^2,\theta^3,\theta^4,\theta^5,\theta^6)$ to
write down the Levi-Civita connection for $\tilde{g}_F$. We have
$$\tilde{g}_F=g_{ij}\theta^i\theta^j,$$
with the indices range: $i,j=0,1,2,...6$, and the matrix $g_{ij}$ given by
$$
g_{ij}=\bma
0&0&0&0&0&0&-1\\
0&0&0&0&0&t^2&0\\
0&0&0&0&-t^2&0&0\\
0&0&0&t^2&0&0&0\\
0&0&-t^2&0&0&0&0\\
0&t^2&0&0&0&-2tu P&0\\
-1&0&0&0&0&0&0
\ema.
$$ 
The Levi-Civita connection for $\tilde{g}_F$ on $J\times\bbR_+\times\bbR$, when
written in the coframe
$(\theta^0,\theta^1,\theta^2,\theta^3,\theta^4,\theta^5,\theta^6)$
reads:
$$
\om_{LC}=\bma
0&0&0&0&0&-tP\theta^5&0\\
&&&&&&\\
\frac{1}{t}\theta^1+\frac{u}{t^2}P\theta^5&\frac{1}{t}\theta^0&Q\theta^2+\tfrac{9}{2\sqrt{3}}P\theta^3&\tfrac{1}{\sqrt{3}}\theta^4&-\tfrac{1}{2\sqrt{3}}\theta^3&\frac{u}{t^2}P\theta^0-\frac{u}{3t}Q\theta^5-\frac{1}{t}P\theta^6&-\frac{1}{t}P\theta^5\\&&&&&&\\
\frac{1}{t}\theta^2&0&\frac{1}{t}\theta^0&\tfrac{1}{\sqrt{3}}\theta^5&0&-\tfrac{1}{2\sqrt{3}}\theta^3&0\\&&&&&&\\
\frac{1}{t}\theta^3&0&-2\sqrt{3}P\theta^5&\frac{1}{t}\theta^0&\tfrac{1}{\sqrt{3}}\theta^5&-\tfrac{1}{\sqrt{3}}\theta^4&0\\&&&&&&\\
\frac1t\theta^4&0&0&-2\sqrt{3}P\theta^5&\frac{1}{t}\theta^0&Q\theta^2+\tfrac{9}{2\sqrt{3}}P\theta^3&0\\&&&&&&\\
\frac{1}{t}\theta^5&0&0&0&0&\frac{1}{t}\theta^0&0\\&&&&&&\\0&t\theta^5&-t\theta^4&t\theta^3&-t\theta^2&t\theta^1-uP\theta^5&0
\ema.
$$ 
Note that on 
$\Sigma=\{(x,y,p,q,z,t,u):u=0,~t=1\}$ we trivially have $\theta^0\equiv
0\equiv\theta^6$.
Thus, restricting the formula for $\om_{LC}$ to $\Sigma$, we see that $\om_{G_2}\equiv \om_{LC|\Sigma}$. Off this
set the two connections: $\om_{LC}$ and the pullbacked-by-$\pi$-connection 
$\om_{G_2}$, differ significantly. To see this it is enough to
observe that contrary to $\om_{LC}$, the connection 
$\pi^*(\om_{G_2})$ has {\it torsion}. Indeed writing the first Cartan structure
equations for the $\pi^*(\om_{G_2})$ in the coframe
$(\theta^0,\theta^1,\theta^2,\theta^3,\theta^4,\theta^5,\theta^6)$ we
find that the torsion is:
$$
\der\theta^i+\pi^*(\om_{G_2})^i_{~j}\dz\theta^j=\bma
0\\
-\theta^0\dz\theta^1-P\theta^5\dz\theta^6\\
-\theta^0\dz\theta^2\\
-\theta^0\dz\theta^3\\
-\theta^0\dz\theta^4\\
-\theta^0\dz\theta^5\\
0
\ema.
$$  
The vanishing of this torsion on the initial hypersurface $\Sigma$
confirms our earlier statemant that the two connections
$\om_{G_2}$ and $\om_{LC}$ coincide there.

It is interesting to note that the curvature
$\der \om_{LC}+\om_{LC}\dz\om_{LC}$ does not depend on $t$, $u$ and is anihilated
by $\partial_t$ and $\partial_u$. Thus it can be considered to be a
2-form on $\Sigma$. As such it is precisely equal to the curvature 
$\der \om_{G_2}+\om_{G_2}\dz\om_{G_2}$ of the connection $\om_{G_2}$:
$$
\der \om_{G_2}+\om_{G_2}\dz\om_{G_2}=\der
\om_{LC}+\om_{LC}\dz\om_{LC}=
\bma
0&0&0&0&0&0&0\\
0&0&A_5&0&0&0&0\\
0&0&0&0&0&0&0\\0&0&0&0&0&0&0\\
0&0&0&0&0&A_5&0\\
0&0&0&0&0&0&0\\
0&0&0&0&0&0&0
\ema\theta^2\dz\theta^5,
$$
where\footnote{We use the letter $A_5$ to denote the nonvanishing
  component of the curvature to be in accordance with \cite{conformal}
and Cartan's paper \cite{cartan}. Note however that in order to avoid
collision of notations between the present and the next sections we use capital $A_5$
instead of $a_5$ of paper \cite{conformal}.}
$$A_5=\frac{-224F^{(3)4}+336F''F^{(3)2}F^{(4)}-51F''^2F^{(4)2}-80F''^2F^{(3)}F^{(5)}+10F''^3F^{(6)}}{100F''^{20/3}}.$$

\section{Non-conformally Einstein example}
To get quite different example of $[g_F]$ we consider equation (\ref{2m}) in the form:
$$z'=y''^2+a_6 y'^6+a_5y'^5+a_4 y'^4+a_3 y'^3+a_2 y'^2+a_1 y'+a_0+b z,$$
where $a_i, i=0,1,...,6,$ and $b$ are real constants. This equation has the defining function
$$F=q^2+a_6 p^6+a_5p^5+a_4 p^4+a_3 p^3+a_2 p^2+a_1 p+a_0+b z$$
and, via (\ref{metr}), leads to a conformal class $[g_F]$ represented
by a metric 
\begin{eqnarray}
&&15 (2)^{-2/3}g_F= [
    9 a_2 + 2b^2+ 
    27 a_3 p + 54 a_4 p^2+ 90 a_5 p^3 + 135 a_6 p^4 ]\der y^2+\nonumber\\&&[15 a_0 + 2 (b^2- 3 a_2) p^2- 3 a_3 p^3 + 9 a_4 p^4+ 30 a_5 p^5+ 60 a_6 p^6 -\nonumber\\&& 20 b  p q + 5 q^2 + 15 b z] \der x^2 +\label{prz2}\\&& [15 a_1 + 4(3 a_2  - b^2) p  - 
    9 a_3 p^2  - 
    48 a_4 p^3  - 
    105a_5 p^4  - 
    180 a_6 p^5 +\nonumber\\&& 
    20 b q]\der x\der y +20\der p^2  -\nonumber
\\&& 10( b p +q)\der p\der x+ 10b \der p \der y - 30\der q\der y   - 15\der x\der z  + 
    30p\der q\der x.\nonumber
\end{eqnarray} 
This metric is \emph{not} conformal to an Einstein metric.  
The quickest way to check this is the calculation of the Cotton,
$C_{ijk}$, and the Weyl, $W_{ijkl}$, tensors for $g_F$. Once these
tensors are calculated, it is easy to observe that they do not admit a
vector field $K^i$ such that $C_{ijk}+K^lW_{lijk}=0$. As a consequence
the metric is \emph{not} a \emph{conformal C-space} metric. This
proves our statement since every conformally Einstein metric is
neccessarily a conformal C-space metric (see e.g. Ref. \cite{gover}). 

Recall that $g_F$ of (\ref{prz2}), as a member  of the family of metrics
(\ref{metr}), defines a conformal class $[g_F]$ with 
\emph{conformal} holonomy $H$ \emph{reduced} to the 
noncompact group $G_2$ or to one of its subgroups. But since the metric
(\ref{prz2}) is 
not conformal to an Einstein metric, 
we do not have an immediate reason to conclude that $H\neq G_2$. 
We \emph{conjecture} that $H=G_2$ here and try to prove it in a
subsequent paper \cite{leis}.

It is remarkable that the ambient metric $\tilde{g}_F$ for $g_F$ of 
(\ref{prz2}) assumes a very compact form:
\begin{eqnarray*}
&&\tilde{g}_F=t^2 g_F~-2~\der t\der u~-~\\
&&2 ~tu~ [~\tfrac{1}{20}(-2 a_2 + 4b^2 + 3 a_3 p + 6 a_4p^2 - 
              20 a_5 p^3 - 120 a_6p^4)\der x^2 -\\&&\tfrac{9}{20}
(a_3 - 10 a_5 p^2 - 40 a_6 p^3)\der x\der y -\tfrac{9}{10}(a_4 + 5 a_5 p + 15a_6 p^2)\der y^2~]~+\\
&&u^2~[~\tfrac{3}{20(2)^{2/3}}(a_4 - 10 a_5 p + 60 a_6 p^2) \der x^2 +
    \tfrac{9}{4(2)^{2/3}}(a_5 - 12a_6 p)\der x\der y +\tfrac{81}{4(2)^{2/3}} a_6 \der y^2~].
\end{eqnarray*}
This is checked by applying our strategy described in
Section \ref{stra} to the metric (\ref{prz2}). As in the previous
example, using Mathematica, we calculated the bilinear form $\gamma$
for (\ref{prz2}). It turned out to be equal to \emph{zero}, $\gamma\equiv
0$. Then we calculated $\bar{g}_F$, and checked that it is
\emph{Ricci flat}. Thus we concluded that
$\bar{g}_F$ coincides with the ambient metric for $\tilde{g}_F$. The above
given formula for $\tilde{g}_F$ is therefore just $\bar{g}_F$,
which we calculated using (\ref{ffgc}).

We find this example as a sort of miracle. Apriori there is no reason
for $g_F$ to have the ambient metric \emph{truncated} at the
\emph{second} order in terms of the ambient parameters $t$ and $u$. We
are intrigued by this fact.

Now, following the general procedure outlined in \cite{conformal}, we
introduce a special coframe for $g_F$ given by:

\begin{eqnarray*}
&&\theta^1=\der y -p\der x\\
&&\theta^2=\der z-F\der x-2q(\der p-q\der x)\\
&&\theta^3=-\tfrac{2^{4/3}}{\sqrt{3}}(\der p-q\der x)\\
&&\theta^4=2^{-1/3}\der x\\
&&15 (2)^{1/3}\theta^5=(9 a_2+2b^2+27 a_3p+54 a_4 p^2+90
  a_5p^3+135a_6p^4)(\der y-p\der x)+\\&&10b(\der p-q\der x)-30\der q+\\
&&15(a_1+2 a_2p+3 a_3 p^2+4 a_4 p^3+5 a_5p^4+6 a_6p^5+2bq)\der x.
\end{eqnarray*}
In this coframe the metric $g_F$ is:
$$g_F=2\theta^1\theta^5-2\theta^2\theta^4+(\theta^3)^2.$$

As in the previous section, we use the canonical projection
$$\pi(x,y,p,q,z,t,u)=(x,y,p,q,z)$$ to pullback 
the coframe $(\theta^1,\theta^2,\theta^3,\theta^4,\theta^5)$ 
to five linearly independent forms
$(\theta^1,\theta^2,\theta^3,\theta^4,\theta^5)$ on
$J\times\bbR_+\times\bbR$, which are further suplemented by $$\theta^0=\der
t\quad\quad{\rm and}\quad\quad \theta^6=\der u$$ to form a coframe
$(\theta^0,\theta^1,\theta^2,\theta^3,\theta^4,\theta^5,\theta^6)$ on
the ambient space $J\times\bbR_+\times\bbR$.

It turns out that if $b= 0$ 
the coframes on $J$ and $J\times\bbR_+\times\bbR$ defined in
this way are suitable to analyze the relations between the Cartan
normal conformal connection $\om_{G_2}$ for $[g_F]$
and the Levi-Civita connection $\om_{LC}$ for $\tilde{g}_F$. If $b\neq
0$
the 
conection $\om_{G_2}$ in the coframe
$(\theta^1,\theta^2,\theta^3,\theta^4,\theta^5)$ and the connection
$\om_{LC}$ in the coframe 
$(\theta^0,\theta^1,\theta^2,\theta^3,\theta^4,\theta^5,\theta^6)$ do
not coincide on $t=1$, $u=0$. We will not
analyze this case here.

Restricting to the 
$$b=0$$
case we find the following:
\begin{itemize}
\item the connections  $\om_{G_2}$ in the coframe
$(\theta^1,\theta^2,\theta^3,\theta^4,\theta^5)$ and the connection
$\om_{LC}$ in the coframe 
$(\theta^0,\theta^1,\theta^2,\theta^3,\theta^4,\theta^5,\theta^6)$ 
coincide on $t=1$, $u=0$.
\item the torsion of $\pi^*(\om_{G_2})$ in the coframe
  $(\theta^0,\theta^1,\theta^2,\theta^3,\theta^4,\theta^5,\theta^6)$
  is nonvanishing off the set $t=1$, $u=0$
\item unlike the example of the previous section the curvature $\der
\om_{LC}+\om_{LC}\dz\om_{LC}$ siginificantly depends on $t$ and $u$. 
\item even on $t=1$, $u=0$, the curvature $\der
  \om_{G_2}+\om_{G_2}\dz\om_{G_2}$ and the restriction of $\der
\om_{LC}+\om_{LC}\dz\om_{LC}$ do not coincide.
\end{itemize}

\section{Acknowledgements}
I am very grateful to T P Branson, M Eastwood and W Miller Jr, the organizers of the 2006 IMA Summer Program
``Symmetries and Overdetermined Systems of Partial Differential
Equations'', for invitating me to Minneapolis to 
particpate in this very fruitful event. The topic covered by this note 
is inspired by the talk of C R Graham which I heard in Minneapolis 
during the program. In particular, I am very obliged to C R Graham for
sending me the formulas (\ref{ffgc}), which I used to prepare the 
examples included in this note.


\begin{thebibliography}{99}
\bibitem{cartan} Cartan E, ``Les systemes de Pfaff a cinq variables et
  les equations aux derivees partielles du seconde ordre''
  \emph{Ann. Sc. Norm. Sup.} {\bf 27} 109-192 (1910)
\bibitem{feff} Fefferman C, Graham C R, ``Conformal invariants'', in
  \emph{Elie Cartan et mathematiques d'aujourd'hui}, Asterisque, hors
  serie (Societe Mathematique de France, Paris) 95-116  (1985)
\bibitem{gover} Gover A R, Nurowski P, "Obstructions to conformally Einstein metrics in n dimensions" {\it Journ. Geom. Phys.}  {\bf 56}  450-484 (2006) 
\bibitem{grah} Graham C R, private communications, unpublished
\bibitem{de} de Haro S, Skenderis K, Solodukhin S N, 
``Holographic reconstruction of spacetime and renormalization in the
AdS/CFT correspondence'', {\it Comm. Math. Phys.} {\bf 217} (2001), 594--622,
hep-th/0002230
\bibitem{leis} Leistner Th, Nurowski P, in preparation
 \bibitem{conformal} Nurowski P, "Differential equations and conformal structures" {\it Journ. Geom. Phys.}  {\bf 55}  19-49 (2005)
\end{thebibliography}
\end{document}